\documentclass[preprint,12pt]{elsarticle}
\usepackage{amsfonts}
\input{epsf.sty}
\usepackage{latexsym,amsmath,amsfonts,amscd, amsthm}
\usepackage{epsfig}
\usepackage{changebar}
\usepackage{pstricks}
\usepackage{pst-plot}
\usepackage{multirow}
\usepackage{subfigure}
\usepackage{placeins}
\usepackage{amssymb}

\usepackage{mathrsfs}
\usepackage{amssymb,amsthm,hyperref,amsmath,graphicx,bm,color,booktabs}
\usepackage{bm}
\usepackage{comment}

\numberwithin{equation}{section} 

\theoremstyle{definition}

\newcommand{\brac}[1]{\left(#1\right)}

\newcommand{\order}{\mathcal{O}}

\newcommand{\ie}{{\it{i.e.}}}
\newcommand{\eg}{{\it{e.g.\ }}}

\newcommand{\dt}{\Delta t}
\newcommand{\dx}{\Delta x}

\allowdisplaybreaks

\journal{Journal of Computational Physics}

\begin{document}

\baselineskip=2pc

\begin{frontmatter}

\title{Boundary treatment of high order Runge-Kutta methods for hyperbolic conservation laws}

\author{Weifeng Zhao}
\ead{wfzhao@ustb.edu.cn}
\address{Department of Applied Mathematics, University of Science and Technology Beijing, Beijing 100083, China}

\author{Juntao Huang\corref{corr}}
\cortext[corr]{Corresponding author}
\ead{huangj75@msu.edu}
\address{Department of Mathematics, Michigan State University, East Lansing, MI 48824, USA}

\author{Steven J. Ruuth}
\ead{sruuth@sfu.ca}
\address{Department of Mathematics, Simon Fraser University, Burnaby, British Columbia, V5A 1S6 Canada.}

\begin{abstract}
In \cite{ZH2019}, we developed a boundary treatment method for
implicit-explicit (IMEX) Runge-Kutta (RK) methods for solving hyperbolic systems with source terms.
Since IMEX RK methods include explicit ones as special cases,  this boundary treatment method naturally applies to explicit methods as well.
In this paper, we examine this boundary treatment method for the case of explicit RK schemes of arbitrary order applied to hyperbolic conservation laws.
We show that the method not only preserves the accuracy of explicit RK schemes but also possesses good stability.
This compares favourably to the inverse Lax-Wendroff method in \cite{TS2010,TWSN2012} where analysis and numerical experiments
have previously verified the presence of order reduction \cite{TS2010,TWSN2012}.
In addition, we demonstrate that our method performs well for strong-stability-preserving (SSP) RK schemes involving negative coefficients and downwind spatial discretizations.
It is numerically shown that when boundary conditions are present and the proposed boundary
treatment is used, that SSP RK schemes with negative coefficients still allow for larger time steps than
 schemes with all nonnegative coefficients.
In this regard, our boundary treatment method is an effective supplement to SSP RK schemes with/without negative coefficients for initial-boundary value problems for hyperbolic conservation laws.
\end{abstract}


\end{frontmatter}

\section{Introduction}
\setcounter{equation}{0}
\setcounter{figure}{0}
\setcounter{table}{0}

To approximate time-dependent partial differential equations (PDEs), it is common to first discretize the spatial derivatives to obtain a large system of time-dependent ordinary differential equations (ODEs). These ODEs are then discretized by suitable time-stepping techniques such as multistep or Runge-Kutta (RK) methods. For hyperbolic conservation laws, strong-stability-preserving (SSP) methods are often applied to maintain desired monotonicity properties of the underlying flow,
particularly in the case of nonsmooth solutions \cite{Shu1988,SO1988,GST2001}. The basic idea of SSP methods is to assume that the forward Euler method is strongly stable under a suitable time stepping  restriction for some norm or semi-norm, and then to construct higher order RK schemes as convex combinations of forward Euler steps with various step sizes. Thanks in part to the provable SSP property, these methods have found substantial use in the time discretization of PDEs.

Many SSP RK schemes are composed exclusively of nonnegative coefficients and the most widely used one is perhaps the three-stage third-order RK method in \cite{Shu1988,SO1988}. However, it is proved in \cite{GS1998} that there are no four-stage fourth-order SSP RK schemes with all nonnegative coefficients. More seriously, it is impossible to construct an explicit SSP RK method of order greater than four with all nonnegative coefficients \cite{RS2002}.  To address this issue, more stages or negative coefficients are taken into consideration. Research along this line can be found in, \eg \cite{GS1998,SR2002,RS2004,R2005,GR2006}. For RK schemes with negative coefficients, the SSP property can also be achieved if the spatial derivatives corresponding to the negative coefficients are approximated by a so-called downwind spatial discretization \cite{Shu1988,GR2006}. Notably, RK schemes involving negative coefficients may allow larger SSP coefficients and consequently can be even more efficient than those with all nonnegative coefficients \cite{GS1998,SR2002,RS2004,R2005,GR2006}.

This paper is concerned with the boundary treatment of general SSP RK schemes with/without negative coefficients for hyperbolic conservation laws.
For SSP RK schemes with negative coefficients,   \cite{Ha_thesis}  applies  the  Navier-Stokes characteristic boundary conditions (NSCBC) method \cite{PL1992}  to the boundary treatment of hyperbolic conservation systems.
For SSP RK schemes with nonnegative coefficients, a popular boundary treatment method is to impose consistent boundary conditions for each intermediate stage \cite{Carpenter1995}. Unfortunately, these boundary conditions are derived only for RK schemes up to third order. This is remedied for fourth-order schemes in \cite{abarbanel1996,pathria1997}, however, the methods therein only apply to one-dimensional (1D) scalar equations or systems with all characteristics flowing into the domain at the boundary. Different from this, we propose in our previous work \cite{ZH2019} to use the RK schemes themselves at the boundary. This idea, combined with an inverse Lax-Wendroff (ILW) procedure in \cite{TS2010,TWSN2012}, preserves the accuracy and good stability of implicit-explicit (IMEX) RK schemes solving hyperbolic systems with source terms. Since IMEX RK methods include explicit ones as special cases, the method in \cite{ZH2019} naturally applies to explicit methods as well.

In this paper, we examine the boundary treatment method in \cite{ZH2019} for explicit high-order SSP RK schemes of hyperbolic conservation laws.
Specifically, we use finite difference WENO schemes on a Cartesian mesh for the spatial discretization, where the corresponding downwind scheme is applied in the case of negative coefficients.
We show that our method applies to general SSP RK schemes with/without negative coefficients.  Furthermore, we show that it preserves the accuracy of the RK schemes and that it has good stability. These nice properties are verified on a selection of linear and nonlinear problems and third- and fourth-order RK schemes with and without negative coefficients. Additionally, when boundary conditions are present and the proposed boundary treatment is used, the SSP RK schemes with negative coefficients still allow larger time steps than those with all nonnegative coefficients. In this regard, the boundary treatment method in \cite{ZH2019} is an effective supplement to the SSP RK schemes with/without negative coefficients for initial-boundary value problems for conservation laws.

In addition, we demonstrate that there exists some order reduction phenomena for the ILW approach in \cite{TS2010,TWSN2012}. In that method \cite{TS2010,TWSN2012}, it is assumed that the numerical solutions in \emph{intermediate} stages all satisfy the PDE and the corresponding consistent boundary conditions \cite{TS2010,TWSN2012}. However, this assumption is expected to result in some error since the solutions at intermediate stages of RK methods are generally not high-order accurate. By contrast, we use the fact that solutions at intermediate stages satisfy a semi-discrete scheme instead of the PDE and then apply the ILW procedure. We analyze the difference between our boundary treatment method and that given in \cite{TS2010,TWSN2012} for the problem of imposing intermediate boundary conditions for the three-stage third-order SSP RK scheme. Our analysis shows the equivalence of the two methods for linear problems. However, for nonlinear equations there exist difference terms of $\order(\dt^2\dx)$ between the two methods. As a consequence, the convergence order is only ${17}/{3}$ for the method in \cite{TS2010,TWSN2012} if using a seventh-order WENO scheme in space and the third-order SSP RK scheme in time with a time step-size $\dt = \order(\dx^{7/3})$. This is less than the expected seventh-order that we obtain with our boundary treatment method.  These interesting findings are verified through numerical experiments as well.

This paper is organized as follows. In Section 2, we introduce the SSP RK schemes and the WENO scheme. In Section 3, we use the one-dimensional case to illustrate our idea of boundary treatment appearing in \cite{ZH2019}. The method is compared with that in \cite{TS2010,TWSN2012} for linear and nonlinear problems in Section 4. Numerical tests are presented in Section 5 to demonstrate the stability and accuracy of  our method as well as the order reduction of the method in \cite{TS2010,TWSN2012}. Finally, Section 6 concludes the paper. This paper also includes an appendix which gives some details of the boundary treatment.

\section{Scheme formulation}
\subsection{RK methods for conservation laws}

Consider a one-dimensional hyperbolic system of conservation laws for $U=U(t,x) \in \mathbb R^M$
\begin{equation}\label{21}
\partial_t U + \partial_x F(U) = 0
\end{equation}
on a bounded domain $0 \leq x \leq 1$ subject to appropriate boundary conditions and initial data.
We assume that the Jacobian matrix $F_{U}(U(t,0))$ always has $p$ positive eigenvalues and thus $p$ independent relations among incoming and outgoing modes are given at the left boundary $x=0$:
\begin{equation}
B(U(t,0), t) = 0.   \label{23a}
\end{equation}
Similarly, the boundary conditions at $x=1$ are also imposed properly.

A general $s$-stage explicit SSP RK scheme for the hyperbolic system \eqref{21} is written as \cite{SO1988,Shu1988}
\begin{subequations}\label{24}
\begin{align}
& U^{(0)} = U^n, \label{24a}\\
& U^{(i)} = \sum_{k=0}^{i-1} \alpha_{ik} U^{(k)}
            + \Delta t \beta_{ik}
            \left\{
               \begin{array}{cc}
                  L( U^{(k)} ) & \mbox{if} ~ \beta_{ik} >0\\
                  \tilde{L}( U^{(k)} )  & \mbox{otherwise}
               \end{array}
            \right. , \quad i=1,2,\ldots, s, \label{24b}\\
& U^{n+1} = U^{(s)}, \label{24c}
\end{align}
\end{subequations}
where all the $\alpha_{ik} \geq 0$ and $\alpha_{ik} = 0$ only if $\beta_{ik} = 0$. Here both $L(U)$ and $\tilde{L}(U)$ approximate $-\partial_x F(U)$, and the forward Euler method
applied to $\dot{U}=L(U)$ is strongly stable under a certain time-step restriction, i.e., $|| U + \Delta t L(U) || \le ||U||$ for all $\Delta t\le \Delta t_{FE}$.
Similarly, the backward-in-time Euler method  applied to $\dot{U}=\tilde{L}(U)$ is strongly stable under a suitable time-step restriction, i.e., $|| U - \Delta t \tilde{L}(U) || \le ||U||$ for all $\Delta t\le \Delta t_{FE}$.
 In practice, $\tilde{L}(U)$ is formed using a reversal in the upwinding direction relative to $L(U)$, or, in other words, using downwinding (see \cite{Shu1988,RS2004,GR2006} for more details). For convenience, we refer to such discretizations as downwind spatial discretizations (cf. \cite{RS2004,GR2006} ).
It is shown in \cite{RS2004,GR2006} that RK methods with negative coefficients may allow larger time steps (or CFL numbers) and even better efficiency than those with all nonnegative coefficients.


\subsection{WENO scheme}
For the spatial discretization, we use the finite difference WENO scheme with the Lax-Friedrichs flux splitting \cite{JS1996}.
Denote by $x_j$ the $j$-th grid point of a uniform mesh and by $U_j$ the corresponding point value. Define $\Delta x = x_{j+1}-x_j$ and $x_{j+\frac{1}{2}} = \frac{1}{2}(x_{j} + x_{j+1})$. The finite difference WENO scheme for the spatial derivative $F(U)_x$ at $x=x_j$
will take the conservative form
\begin{equation*}
 \frac{1}{\Delta x}( \hat F_{j+\frac{1}{2}} - \hat F_{j-\frac{1}{2}} ) ,
\end{equation*}
where $\hat F_{j+\frac{1}{2}}$ is the numerical flux. Namely, $L(U)$ in \eqref{24} for $\beta_{ik}>0$ is given by
\begin{equation}
L(U) = - \frac{1}{\Delta x}( \hat F_{j+\frac{1}{2}} - \hat F_{j-\frac{1}{2}} ).
\end{equation}
We use the global Lax-Friedrichs splitting
$$
F_{\pm}(U) := \frac{1}{2} \left( U \pm \frac{F(U)}{\alpha} \right),
$$
where $\alpha = \max{ |\lambda_i |}$ with $\lambda_i$ being the $i$-th eigenvalue of the Jacobian matrix $F_{U}(U)$.  The maximum is taken over
all the grid points at time level $t_n$. Let $L_{j+\frac{1}{2}}$ and $R_{j+\frac{1}{2}}$ be the left and right eigenvector matrix of $F_U(U)$ evaluated at $U_{j+\frac{1}{2}}$.  These satisfy $L_{j+\frac{1}{2}}^{-1} = R_{j+\frac{1}{2}}$. Here $U_{j+\frac{1}{2}}$ is some average of $U_j$
and $U_{j+1}$, and we simply take $U_{j+\frac{1}{2}}=\frac{1}{2}( U_j + U_{j+1} )$. The finite difference WENO scheme is formulated as follows.

At each fixed $x_{j+\frac{1}{2}}$:

1. Transform the cell average $F_{\pm}(U_k^n)$ for all $k$ in a neighborhood of $j$ to the local characteristic field by setting
$$
(V_{\pm})_{k}^n := L_{j+\frac{1}{2}} F_{\pm}(U_k^n).
$$

2. Perform the WENO reconstruction for each component of $(V_{+})_{k}^n$ to obtain the value of $(V_{+})^n$ at $x_{j+\frac{1}{2}}$ on the left side
and denote it as $(V_{+})^{-}_{j+\frac{1}{2}}$. Similarly, perform the WENO reconstruction for each component of $(V_{-})_{k}^n$ to obtain the value of $(V_{-})^n$ at $x_{j+\frac{1}{2}}$ on the right side and denote it as $(V_{-})^{+}_{j+\frac{1}{2}}$.

3. Transform back to the physical space by
$$
(F_{+})^{-}_{j+\frac{1}{2}} := R_{j+\frac{1}{2}} (V_{+})^{-}_{j+\frac{1}{2}}, \qquad
(F_{-})^{+}_{j+\frac{1}{2}} := R_{j+\frac{1}{2}} (V_{-})^{+}_{j+\frac{1}{2}}.
$$

4. Form the flux using
$$
\hat F_{j+\frac{1}{2}} = \alpha\left[  (F_{+})^{-}_{j+\frac{1}{2}} - (F_{-})^{+}_{j+\frac{1}{2}} \right].
$$

\subsection{Downwind space discretization}

To implement the downwind spatial discretization for $\tilde{L}(U)$ in the case of $\beta_{ik}<0$, we only need to define \cite{GR2006}
\begin{equation}
\tilde{F}_{+}(U) := F_{-}(U) = \frac{1}{2} \left( U - \frac{F(U)}{\alpha} \right) , \quad
\tilde{F}_{-}(U) := F_{+}(U) = \frac{1}{2} \left( U + \frac{F(U)}{\alpha} \right).
\end{equation}
Using these expressions,  the procedure to compute the flux is the same as that of the WENO scheme above. For convenience, we provide the details as well. At each fixed $x_{j+\frac{1}{2}}$:

1. Transform the cell average $\tilde{F}_{\pm}(U_k^n)$ for all $k$ in a neighborhood of $j$ to the local characteristic field by setting
$$
(\tilde{V}_{\pm})_{k}^n := L_{j+\frac{1}{2}} \tilde{F}_{\pm}(U_k^n).
$$

2. Perform the WENO reconstruction for each component of $(\tilde{V}_{+})_{k}^n$ to obtain the value of $(\tilde{V}_{+})^n$ at $x_{j+\frac{1}{2}}$ on the left side and denote it as $(\tilde{V}_{+})^{-}_{j+\frac{1}{2}}$. Similarly, perform the WENO reconstruction for each component of $(\tilde{V}_{-})_{k}^n$ to obtain the value of $(\tilde{V}_{-})^n$ at $x_{j+\frac{1}{2}}$ on the right side and denote it as $(\tilde{V}_{-})^{+}_{j+\frac{1}{2}}$.

3. Transform back to the physical space by
$$
(\tilde{F}_{+})^{-}_{j+\frac{1}{2}} := R_{j+\frac{1}{2}} (\tilde{V}_{+})^{-}_{j+\frac{1}{2}}, \qquad
(\tilde{F}_{-})^{+}_{j+\frac{1}{2}} := R_{j+\frac{1}{2}} (\tilde{V}_{-})^{+}_{j+\frac{1}{2}}.
$$

4. Form the flux using
$$
\tilde{\hat F}_{j+\frac{1}{2}} = -\alpha\left[  (\tilde{F}_{+})^{-}_{j+\frac{1}{2}} - (\tilde{F}_{-})^{+}_{j+\frac{1}{2}} \right].
\quad (\mbox{note the additional minus sign})
$$

With the above flux $\tilde{\hat F}_{j+\frac{1}{2}}$, $\tilde{L}(U)$ in \eqref{24} for $\beta_{ik}<0$ is computed as
\begin{equation}
\tilde{L}(U) =  -\frac{1}{\Delta x}( \tilde{\hat F}_{j+\frac{1}{2}} - \tilde{\hat F}_{j-\frac{1}{2}} ).
\end{equation}

\section{Boundary treatment}
\label{sec3}

In this section, we introduce the boundary treatment method in \cite{ZH2019} with the fifth-order finite difference WENO scheme. Here we take $\Delta t=\order(\Delta x)$.

\subsection{Computation of solutions at ghost points}\label{sec31}

We focus on the left boundary $x=0$ of the problem \eqref{21}; the method can be similarly applied to the right boundary.
Following the notations in \cite{TWSN2012}, we discretize the interval $[0,1]$ by a uniform mesh
\begin{equation}\label{3-1}
\frac{\Delta x}{2} = x_0 < x_1 < \cdots <x_N = 1- \frac{\Delta x}{2}
\end{equation}
and set $x_j, j=-1,-2,-3,$ as three ghost points near the left boundary $x=0$ (note that the boundary can be arbitrarily located).
Denote by $U_j^n$ the numerical solution of $U$ at position $x_j$ and time $t_n$.
Assume that the interior solutions $U_j, j=0,1,2,\ldots,N$, have been updated from time level $t_{n-1}$ to time level $t_n$.
Since the spatial discretization is fifth-order, we use a fifth-order Taylor approximation to construct the values at the ghost points,
\begin{equation}\label{25}
U_j^n = \sum_{k=0}^{4} \frac{x_j^k}{k!} U^{n,(k)} , \quad j=-1,-2,-3,
\end{equation}
where $U^{n,(k)}$ denotes a $(5-k)$-th order approximation of the spatial derivative at the boundary point $\frac{\partial^k U}{\partial x^k}\big|_{x=0, t=t_n}$.
With this formula, $U_j^n$ at the ghost points can be obtained once $U^{n,(k)}, k=0,1,2,3,4$ are provided.  Full details on this procedure appear in \cite{TS2010,TWSN2012} (see also the next subsection).
Having found the ghost values, we can obtain  $U_j^{(1)}$ at the interior points via the RK scheme \eqref{24}.

The next task is to compute $U_j^{(1)}, j=-1,-2,-3$, which is a key point of our method.
Similar to the above procedure, we compute $U_j^{(1)}$ using the fifth-order Taylor expansion at the boundary point $x_b=0$:
\begin{equation}\label{27}
U_j^{(1)} = \sum_{k=0}^{4} \frac{x_j^k}{k!} U^{(1),(k)},  \quad j=-1,-2,-3,
\end{equation}
where $U^{(1),(k)}$ denotes a $(5-k)$-th order approximation of the spatial derivative $\frac{\partial^k U^{(1)}}{\partial x^k}\big|_{x=0}$.
Next, we compute $U^{(1),(k)}$ for $k=0,1,2,3,4$.
To this end, we apply the first stage of the RK solver \eqref{24b} for $U^{(1)}$ at the boundary point $x_b=0$:
\begin{equation}\label{28}
U^{(1)}(x_b) = \alpha_{10} U^{n}(x_b) - \beta_{10}  \Delta t \partial_x F(U^n(x_b)).
\end{equation}
Notice that $U^{n}(x_b) = U^{n,(0)} + \order(\Delta x^5)$ and
$\partial_x F(U^n(x_b)) = F_U(U^n(x_b)) \partial_x U^n(x_b) = F_U(U^{n,(0)}) U^{n,(1)} + \order(\Delta x^4)$ are already known.
Substituting these into \eqref{28}, we obtain an approximation of  $U^{(1)}(x_b)$ and denote it by $U^{(1),(0)}$. Observe that the error between $U^{(1),(0)}$ and  $U^{(1)}(x_b)$ defined by \eqref{28} is $\order(\Delta x^5)$.


Taking derivatives with respect to $x$ on both sides of \eqref{28} yields
\begin{equation}\label{29}
\frac{\partial U^{(1)}}{\partial x}\big|_{x= x_b}
= \alpha_{10} \partial_x U^{n} (x_b)
 -\beta_{10} \Delta t \partial_{xx} F(U^n(x_b)).
\end{equation}
Here $\partial_x U^{n} (x_b) = U^{n,(1)} + \order(\Delta x^4)$ and
$\partial_{xx} F(U^n(x_b)) = F_{UU}(U^n(x_b)) \partial_x U^n(x_b) \partial_x U^n(x_b) + F_U(U^n(x_b)) \partial_{xx} U^n(x_b)
= F_{UU}(U^{n,(0)}) U^{n,(1)} U^{n,(1)} + F_U(U^{n,(0)}) U^{n,(2)} +  \order(\Delta x^3)$.
With these approximations, $\frac{\partial U^{(1)}}{\partial x}\big|_{x= x_b} $ can be computed with \eqref{29} and the resulting solution  $U^{(1),(1)}$ is a fourth-order approximation of $\frac{\partial U^{(1)}}{\partial x}\big|_{x= x_b} $.

By taking higher order derivatives on both sides of \eqref{28}, one can also compute $\frac{\partial ^k U^{(1)}}{\partial x^k}\big|_{x= x_b}$ for $k=2$. However, this procedure is quite complicated as it involves the Jacobian of a Jacobian. Here we simply approximate $\frac{\partial ^k U^{(1)}}{\partial x^k}\big|_{x= x_b} $ for $k\geq2$ by using a $(5-k)$-th order extrapolation with $U^{(1)}$ at interior points.
In this way, we obtain $U^{(1),(k)}$ for $k=0,1,2,3,4$ with accuracy of order $(5-k)$. Then, $U_j^{(1)}$ for $j=-1,-2,-3$ can be computed by \eqref{27} with fifth-order accuracy. Having $U^{(1)}$ at the ghost points, we can then evolve from $U^{(1)}$ to $U^{(2)}$ using the interior difference scheme.

Repeating the same procedure for each $U^{(i)}, 2\le i\le s-1$, we can compute the solution at the ghost points in the $i$-th intermediate stage. Following this, we can update the solution at all interior points in the $(i+1)$-th stage. Finally, we obtain $U^{n+1}=U^{(s)}$, \ie, the solution at the end of a complete RK cycle. { For clarity, we provide the computation of $U_j^{(2)}$ at ghost points $j=-1,-2,-3$ in the Appendix.}

\subsection{Computation of $U^{n,(k)}$ at the boundary}\label{sec32}

For the sake of completeness, we provide the method in \cite{TWSN2012} for computing $U^{n,(k)}, 0\le k\le 4$, \ie, the $(5-k)$-th order approximation of $\frac{\partial^k U^n}{\partial x^k}$ at the boundary $x=0$.

\subsubsection{$k=0$}


We first do a local characteristic decomposition to determine the inflow and outflow boundary conditions as in \cite{TWSN2012}. Denote the Jacobian matrix of the flux evaluated at $x=x_0$ by
$$A(U_0^n) = \partial_U F(U)\big|_{U=U_0^n}$$
and assume that it has $p$ positive eigenvalues $\lambda_1, \lambda_2, \ldots, \lambda_p$ and $(M-p)$ negative eigenvalues $\lambda_{p+1}, \lambda_{p+2}, \ldots, \lambda_M$ with $l_1, l_2, ...l_p$ and $l_{p+1}, l_{p+2}, ...l_M$ the corresponding left eigenvectors, respectively.
Define by $V_{j,m}$ the $m$-th component of the local characteristic variable at grid point $x_j,j=0,1,2,3,4$, \ie,
\begin{equation}\label{310}
V_{j,m} = l_m U_j^n, \quad  m=1,2,\ldots,M, \quad j=0,1,2,3,4.
\end{equation}
We extrapolate the outgoing characteristic variable $V_{j,m}, m=p+1,p+2,\ldots,M$ to the boundary $x_b$ (for smooth solutions, we use Lagrangian extrapolation, otherwise the WENO type extrapolation in \cite{TS2010,TWSN2012} can be used to avoid possible oscillations), and denote the extrapolated $k$-th order derivative by
\begin{equation}\label{311}
V_{x_b, m}^{*(k)}, \quad k=0,1,2,3,4.
\end{equation}
With $V_{x_b, m}^{*(0)}$ and the boundary condition \eqref{23a}, $U^{n,(0)}$ at the boundary can be determined from the following equations
\begin{equation}\label{312_0}
\begin{split}
& l_m U^{n,(0)} = V_{x_b, m}^{*(0)},  \quad m = p+1,p+2,\ldots,M, \\
& B( U^{n,(0)}, t_n ) = 0.
\end{split}
\end{equation}

\subsubsection{$k=1$}\label{sec322}

Having $U^{n,(0)}$, we proceed to compute $U^{n,(1)}$ with the ILW procedure proposed in \cite{TS2010,TWSN2012}.
To do this, we differentiate the boundary condition \eqref{23a} with respect to $t$
\begin{equation*}
B_U( U^{n,(0)}, t_n  ) \partial_t U^{n,(0)} + B_t( U^{n,(0)}, t_n  ) = 0,
\end{equation*}
which can be written as
\begin{equation*}
B_U( U^{n,(0)}, t_n  ) \partial_t U^{n,(0)} = g(U^{n,(0)} , t_n )
\end{equation*}
with $g(U^{n,(0)} , t_n ) := -B_t( U^{n,(0)}, t_n  )$. In addition, multiplying the equation \eqref{21} with $B_U( U^{n,(0)}, t_n  )$ from the left yields
$$
B_U( U^{n,(0)}, t_n  ) \partial_t U^{n,(0)}  +  B_U( U^{n,(0)}, t_n  ) A( U^{n,(0)} ) U^{n,(1)} = 0.
$$
With the above two equations and $V_{x_b, m}^{*(1)}$ obtained by the extrapolation, $U^{n,(1)}$ can be determined by solving
\begin{subequations}\label{312}
\begin{align}
& l_m U^{n,(1)} = V_{x_b, m}^{*(1)},  \quad m = p+1,p+2,\ldots,M,  \label{312a}\\
& B_U( U^{n,(0)}, t_n  ) A( U^{n,(0)} ) U^{n,(1)} =  - g(U^{n,(0)} , t_n ). \label{312b}
\end{align}
\end{subequations}
This is the ILW procedure in \cite{TS2010,TWSN2012}.


\subsubsection{$k\geq2$}
Following \cite{TWSN2012}, we simply extrapolate $\frac{\partial^k U^n}{\partial x^k}|_{x=0}$ for all $k\geq2$, a procedure that will not affect the stability \cite{TWSN2012}. First, note that we have already obtained the characteristic variables at grid points near the boundary in \eqref{310}.  Based on these characteristic variables, Lagrangian extrapolation (for smooth solutions) or the WENO type extrapolation (for nonsmooth solutions)  is employed to compute the derivatives $V_{x_b, m}^{*(k)}, k\geq2$ for each $m=1,2,\cdots, M$. Following this step, the approximation of $\frac{\partial^k U^n}{\partial x^k}|_{x=0}$ is given by
\begin{equation}\label{3}
U^{n,(k)} = l^{-1} V_{x_b, m}^{*(k)},
\end{equation}
where $l$ is the local left eigenvector matrix composed of $l_1, l_2, \cdots, l_M$.

Our boundary treatment for the general RK scheme \eqref{24} is summarized as follows:

{\bf Step 1:} Compute $U^{n,(k)}$, the $(5-k)$-th order approximation of $\frac{\partial^k U^n}{\partial x^k}$ at the boundary $x=0$ for $k=0,1,2,3,4$ as in subsection \ref{sec32}. Next, impose grid values $U^n_j$ at the ghost points $j=-1,-2,-3$ using the Taylor expansion \eqref{25}. This step is the same as that in \cite{TWSN2012}. With $U^n$ prescribed at the ghost points, we can evolve from $U^n$ to $U^{(1)}$ with the interior difference scheme.

{\bf Step 2:} Compute $U^{(1),(k)}$, the $(5-k)$-th order approximation of $\frac{\partial^k U^{(1)}}{\partial x^k}$ at the boundary $x=0$ for $k=0,1,2,3,4$, as in subsection \ref{sec31}. Specifically, $U^{(1),(0)}$ and $U^{(1),(1)}$ are obtained from the RK solver via \eqref{28} and \eqref{29}, respectively, and $U^{(1),(k)}, k=2,3,4$, are simply extrapolated. Next, impose grid values $U^{(1)}_j$ at the ghost points $j=-1,-2,-3$ with the Taylor expansion \eqref{27}. Having $U^{(1)}$ at the ghost points, we can update $U^{(2)}$ at all interior points.

{\bf Step 3:} Repeat the same procedure as in Step 2 for each $U^{(i)}, 2\le i\le s-1,$ to compute the solution at the ghost points in the $i$-th intermediate stage. We can then update the solution at all interior points in the $(i+1)$-th stage. Finally, we obtain $U^{n+1}=U^{(s)}$, \ie, the solution at the end of a complete RK cycle.

\section{Comparison with the method in \cite{TWSN2012}}

In this section, we analyze the difference between 
our method and that in \cite{TWSN2012}, where the intermediate boundary conditions in \cite{Carpenter1995} are employed.
To this end, we consider the scalar conservation law
\begin{equation}\label{scalar}
	u_t + f(u)_x = 0, \quad 0<x<1
\end{equation}
with boundary condition
\begin{equation}
	u(0,t)=g(t).
\end{equation}
Here $f'(u(0,t))>0$ is assumed. We solve the problem with the third-order SSP RK scheme \cite{SO1988,Shu1988}
\begin{subequations}\label{SSP3}
\begin{align}
& u^{(1)} = u^n + \dt \mathcal{L}( u^n  ), \label{4-00a} \\
& u^{(2)} = \frac{3}{4}u^n + \frac{1}{4}u^{(1)}
               + \frac{1}{4} \dt \mathcal{L}( u^{(1)}  ), \label{4-00b}\\
& u^{n+1} = \frac{1}{3}u^n + \frac{2}{3}u^{(2)}
               + \frac{2}{3} \dt \mathcal{L}( u^{(2)}  ), \label{4-00c}
\end{align}
\end{subequations}
where $\mathcal{L}( u ) = - f(u)_x$.

The method in \cite{TWSN2012} assumes that $u^n$, $u^{(1)}$ and $u^{(2)}$ all satisfy the PDE \eqref{scalar} together with the boundary conditions \cite{Carpenter1995}
\begin{subequations}\label{SSP3-bc}
\begin{align}
	& u^n = g(t^n), \\
	& u^{(1)} = g(t^n) + \dt g'(t^n), \label{44b}\\
	& u^{(2)} = g(t^n) + \frac{1}{2} \dt g'(t^n) + \frac{1}{4} (\dt)^2 g''(t^n).
\end{align}	
\end{subequations}
Then, the first-order derivatives $u^n_x$, $u^{(1)}_x$ and $u^{(2)}_x$ at the boundary $x=0$ are obtained by the ILW procedure:
\begin{equation}\label{45}
	u^n_x = -\frac{u^n_t}{f'(u^n)} = -\frac{g'(t^n)}{f'(g(t^n))},
\end{equation}
\begin{equation}\label{46}
	u^{(1)}_x = -\frac{u^{(1)}_t}{f'(u^{(1)})} = -\frac{(g(t^n) + \dt g'(t^n))'}{f'(g(t^n) + \dt g'(t^n))} = -\frac{g'(t^n) + \dt g''(t^n)}{f'(g(t^n) + \dt g'(t^n))},
\end{equation}
\begin{equation}\label{47}
	u^{(2)}_x = -\frac{u^{(2)}_t}{f'(u^{(2)})} = -\frac{(g(t^n) + \frac{1}{2} \dt g'(t^n) + \frac{1}{4} (\dt)^2 g''(t^n))'}{f'(g(t^n) + \frac{1}{2} \dt g'(t^n) + \frac{1}{4} (\dt)^2 g''(t^n))} = -\frac{g'(t^n) + \frac{1}{2} \dt g''(t^n) + \frac{1}{4} (\dt)^2 g'''(t^n)}{f'(g(t^n) + \frac{1}{2} \dt g'(t^n) + \frac{1}{4} (\dt)^2 g''(t^n))}.
\end{equation}

%


\subsection{Linear case}
We first consider the linear case with $f(u)=u$ .  Here, the derivatives \eqref{45}--\eqref{47} for the method in \cite{TWSN2012} reduce to
\begin{subequations}
\begin{align}
	& u_x^n = -g'(t^n), \label{48a} \\
	& u_x^{(1)} = -g'(t^n) - \dt g''(t^n), \label{48b} \\
	& u_x^{(2)} = -g'(t^n) - \frac{1}{2} \dt g''(t^n) - \frac{1}{4} (\dt)^2 g^{'''}(t^n). \label{48c}
\end{align}	
\end{subequations}
Instead of imposing intermediate boundary conditions, our method directly uses the RK scheme \eqref{SSP3} at the boundary. This means
\begin{subequations}
	\begin{align*}
		u_x^{(1)} ={}& u_x^n + \dt u_{tx}^n, \\
				={}& - u_t^n - \dt u_{tt}^n, \\
				={}& -g'(t^n) - \dt g''(t^n),
	\end{align*}
\end{subequations}
which is the same as \eqref{48b}. Taking one more derivative gives
\begin{subequations}
	\begin{align*}
		u_{xx}^{(1)} ={}& u_{xx}^n + \dt u_{txx}^n, \\
				={}& u_{tt}^n + \dt u_{ttt}^n, \\
				={}& g''(t^n) + \dt g^{'''}(t^n)
	\end{align*}
\end{subequations}
and thereby
\begin{subequations}
	\begin{align*}
		u_x^{(2)} ={}& \frac{3}{4}u_x^n + \frac{1}{4}u_x^{(1)} - \frac{1}{4}\dt u_{xx}^{(1)}, \\
					={}& -\frac{3}{4}g'(t^n) + \frac{1}{4}(-g'(t^n) - \dt g''(t^n)) - \frac{1}{4}\dt (g''(t^n) + \dt g^{'''}(t^n)), \\
					={}& -g'(t^n) - \frac{1}{2} \dt g''(t^n) - \frac{1}{4}\dt^2 g^{'''}(t^n),
	\end{align*}
\end{subequations}
which is the same as \eqref{48c}. Thus, our method and that in \cite{TWSN2012} are equivalent for linear equations.

\subsection{Nonlinear case}
Next, we consider the nonlinear case. Since our boundary treatment of $u^n$ is the same as that in \cite{TWSN2012}, we have
\begin{equation}\label{412}
	u^n_x = -\frac{g'(t^n)}{f'(g(t^n))},
\end{equation}
\begin{equation}
u^n_{xt} = -\frac{u^n_{tt}f'(u^n)-(u^n_t)^2f''(u^n)}{(f'(u^n))^2}
= -\frac{g''(t_n)f'(g(t_n))-(g'(t_n))^2f''(g(t_n))}{(f'(g(t_n)))^2}.
\end{equation}
%
For $u^{(1)}$, the value at the boundary in our method is determined by the RK scheme \eqref{4-00a}  (instead of \eqref{44b}).  Thus
\begin{equation}\label{414}
	u^{(1)} = u^{(n)} + \dt u^{(n)}_t = g(t^n)+\dt g'(t^n).
\end{equation}
It follows from \eqref{412}--\eqref{414} that the first order derivative $u^{(1)}_x$ is
\begin{subequations}
	\begin{align*}
		u_x^{(1)} ={}& u_x^n + \dt u_{xt}^n, \\
				={}& -\frac{u_t^n}{f'(u^n)} + \dt\frac{-u_{tt}^nf'(u^n)+(u_t^n)^2f''(u^n)}{(f'(u^n))^2}, \\
				={}& -\frac{g'(t^n)f'(g(t^n)) + \dt g''(t^n)f'(g(t^n)) - \dt(g'(t^n))^2 f''(g(t^n)) }{f'(g(t^n))^2}.
	\end{align*}	
\end{subequations}
We make a subtraction between the above $u_x^{(1)}$ and that in \eqref{46}:
\begin{subequations}
	\begin{align*}
		{}& -\frac{g'(t) + \dt g''(t)}{f'(g(t) + \dt g'(t))} + \frac{g'(t)f'(g(t)) + \dt g''(t)f'(g(t)) - \dt(g'(t))^2 f''(g(t)) }{f'(g(t))^2} \\
		={}& \frac{2f'(g(t)) f''(g(t)) g'(t) g''(t) + (g'(t))^3\brac{-2(f''(g(t)))^2 + f'(g(t))f^{(3)}(g(t))}}{2(f'(g(t)))^3}\dt^2 + \order(\dt^3).
	\end{align*}
\end{subequations}
It is clear that there exists a difference of $\order(\dt^2)$ between the values of $u_x^{(1)}$ used by the two methods.
This implies that the solutions $u^{(1)}$ at ghost points computed by Taylor expansion will differ by $\order(\dt^2\dx)$.
Similarly, we can show that the difference between values of $u_x^{(2)}$ is also $\order(\dt^2)$ and that the solutions $u^{(2)}$ at ghost points will differ by $\order(\dt^2\dx)$ as well.

The differences computed above show that if we use a seventh-order spatial discretization and take $\dt=\order(\dx^{7/3})$, then values of $u^{(1)}$ (and also $u^{(2)}$) at ghost points will differ by $\order(\dx^{17/3})$
between the two methods. Consequently, order reduction will arise for the method in \cite{TWSN2012}:  the convergence order is only ${17}/{3}$ which is less than the expected seventh-order.
Numerical verification of these results appears in the next section.



\section{Numerical validations}

This section is devoted to the validation of our boundary treatment method.  Unless otherwise stated, we use the fifth-order WENO scheme \cite{SO1988} for the spatial discretization. Fifth-order Lagrangian extrapolation is employed except for the blast wave problem of the Euler equations.  For time-stepping, four RK schemes are considered: the three-stage third-order SSP RK scheme in \cite{SO1988}, the optimal five-stage fourth-order RK scheme with all nonnegative coefficients \cite{R2005}, and the three-stage third-order and five-stage fourth-order RK schemes with negative coefficients proposed in \cite{GR2006}. These four schemes are referred to as SSP(3,3), SSP(5,4), $\mbox{SSP}_{*}(3,3)$ and $\mbox{SSP}_{*}(5,4)$, respectively \cite{R2005,GR2006}.
The coefficients defining the schemes and the corresponding SSP coefficients are given in Table \ref{tb1}.
It should be noted that $\mbox{SSP}_{*}(3,3)$  and $\mbox{SSP}_{*}(5,4)$ are optimal under the restriction that both
$L(U)$ and $\tilde{L}(U)$ arise simultaneously at one level $j$ only, \ie, $\beta_{ij}, j+1 \le i \le s,$ can be positive or negative (or zero)
while $\beta_{ik} \geq 0$ for $k\neq j$.
\begin{table}[!htbp]\centering
\caption{The coefficients of four RK schemes.}
{\tiny
\begin{tabular}{llc}
\toprule
scheme    &  coefficients & SSP coefficient \\ \hline
SSP(3,3)  &
           $
           \begin{aligned}
            & \alpha_{ik} &
                 \begin{array}{ccc}
                  1 & & \\
                  \frac{3}{4}  & \frac{1}{4} & \\
                  \frac{1}{3} & 0 & \frac{2}{3}
                 \end{array}
            \\ \hline
            & \beta_{ik} &
                 \begin{array}{ccc}
                  1 & & \\
                  0 &  \frac{1}{4} & \\
                  0 & 0 & \frac{2}{3}
                 \end{array}
            \end{aligned}
           $
         & 1          \\ \hline
$\mbox{SSP}_{*}(3,3)$
         & $
           \begin{aligned}
            & \alpha_{ik} &
                 \begin{array}{ccc}
                  1 & & \\
                  0.410802706918667  & 0.589197293081333 & \\
                  0.123062611901395 & 0.251481201947289 & 0.625456186151316
                \end{array}
            \\ \hline
           & \beta_{ik}  &
                 \begin{array}{ccc}
                  0.767591879243998 & & \\
                  -0.315328821802221 &  0.452263057441777 & \\
                  -0.041647109531262 &  0 & 0.480095089312672
                 \end{array}
           \end{aligned}
           $
          &   1.3027756          \\ \hline
SSP(5,4)
       & $
           \begin{aligned}
            & \alpha_{ik} &
                 \begin{array}{ccccc}
                  1 & & & &\\
                  0.444370493651235  & 0.555629506348765 & & & \\
                  0.620101851488403 & 0 & 0.379898148511597 & & \\
                  0.178079954393132 & 0 & 0 & 0.821920045606868 & \\
                  0 & 0 & 0.517231671970585 & 0.096059710526147 & 0.386708617503269
                \end{array}
            \\ \hline
           & \beta_{ik}  &
                 \begin{array}{ccccc}
                  0.391752226571890 & & & &\\
                  0 & 0.368410593050371 & & &\\
                  0 & 0 & 0.251891774271694 & &\\
                  0 & 0 & 0 & 0.544974750228521 & \\
                  0 & 0 & 0 & 0.063692468666290 & 0.226007483236906
                 \end{array}
           \end{aligned}
           $
          &   1.5081800         \\ \hline
$\mbox{SSP}_{*}(5,4)$
         & $
           \begin{aligned}
            & \alpha_{ik} &
                 \begin{array}{ccccc}
                  1 & & & &\\
                  0.210186660827794  & 0.789813339172206 & & & \\
                  0.331062996240662 & 0.202036516631465 & 0.466900487127873 & & \\
                  0 & 0 & 0 & 1 & \\
                  0.097315407775058 & 0.435703937692290 & 0 & 0 & 0.466980654532652
                \end{array}
            \\ \hline
           & \beta_{ik}  &
                 \begin{array}{ccccc}
                  0.416596471458169 & & & &\\
                  -0.103478898431154 & 0.388840157514713 & & &\\
                  -0.162988621767813 & 0 & 0.229864007043460 & &\\
                  0 & 0 & 0 & 0.492319055945867 & \\
                  -0.047910229684804 & 0.202097732052527 & 0 & 0 & 0.229903474984498
                 \end{array}
           \end{aligned}
           $
          &   2.0312031         \\
\bottomrule
\end{tabular}\label{tb1}}
\end{table}

With these four RK schemes, we approximate three 1D problems from \cite{TS2010} to test the accuracy and stability of our boundary treatment method.
The order reduction of the method in \cite{TWSN2012} is verified as well.

\subsection{1D linear scalar equation}

The first problem we consider is the linear scalar hyperbolic equation \cite{TS2010} where a boundary condition is prescribed at $x=-1$:
\begin{equation}\label{51}
\begin{aligned}
& u_t + u_x = 0,                     &\quad x \in [-1,1], t>0,\\
& u(0, x) = 0.25 + 0.5 \sin (\pi x), &\quad x\in[-1, 1], \\
& u(t, -1) = g(t),                   &\quad t>0,
\end{aligned}
\end{equation}
No boundary condition is needed at $x=1$.
As in \cite{TS2010}, we first take
\begin{equation}
g(t) = 0.25 - 0.5\sin( \pi(1+t) )
\end{equation}
so that the initial-boundary value problem \eqref{51} has the exact solution
\begin{equation}\label{53}
u(t, x) = 0.25 + 0.5\sin( \pi(x-t) ).
\end{equation}
With this analytical solution, we conduct a numerical convergence study for our boundary treatment method with different RK schemes.
Using a fixed CFL number of 0.6, the $L^1$ and $L^{\infty}$ errors at $t=1$ are computed; see Table \ref{tb52}.
We observe that the desired third-order convergence is obtained for both the SSP(3,3) and $\mbox{SSP}_{*}(3,3)$ schemes. For the fourth-order RK schemes, the convergence orders are about five, an effect which may be due to the error of the spatial discretization dominating over that of the time discretization.

Next, we use the seventh-order WENO scheme and the SSP(3,3) scheme with $\dt = \dx^{7/3}$ to investigate the accuracy of our method and that in \cite{TWSN2012}. As shown in  Table \ref{tb53}, the errors of the two methods are very close to each other for different $\dx$ and the convergence order is about seven for both methods. These results agree with our earlier analysis showing that the method in \cite{TWSN2012} is equivalent to ours for linear equations. It should be noted that we use the ideal weights, which are also called the linear weights, for the seventh-order WENO scheme; using the original ones in \cite{BS2000}, it seems difficult to achieve the intended seventh-order accuracy in practice \cite{SZ2008}.

\begin{table}[!htbp]\centering
\caption{Error table for the linear equation \eqref{51} with solution \eqref{53}. }
\begin{tabular}{lllll | llll}
\toprule
         &  SSP(3,3)  &      &      &      &   $\mbox{SSP}_{*}(3,3)$ &   &  \\\hline
$\Delta x$ & $L^1$ error & order   & $L^{\infty}$ error   & order & $L^1$ error & order  & $L^{\infty}$ error   & order \\ \hline
1/20      &  3.45e-5    &           & 7.44e-5   &        & 3.83e-5     &         & 8.24e-5     &      \\
1/40      &  3.51e-6    & 3.30      & 7.39e-6   &  3.33  & 3.63e-6     & 3.40    & 7.62e-6     & 3.43\\
1/80      &  4.16e-7    & 3.08      & 8.71e-7   &  3.08  & 4.20e-7     & 3.11    & 8.78e-7     & 3.11\\
1/160     &  5.12e-8    & 3.02      & 1.07e-7   &  3.03  & 5.14e-8     & 3.03    & 1.07e-7     & 3.04\\
1/320     &  6.39e-9    & 3.00      & 1.34e-8   &  3.00  & 6.39e-9     & 3.00    & 1.34e-8     & 3.00\\
\toprule
         &  SSP(5,4)  &      &      &      &   $\mbox{SSP}_{*}(5,4)$ &   &  \\\hline
$\Delta x$ & $L^1$ error & order   & $L^{\infty}$ error   & order & $L^1$ error & order  & $L^{\infty}$ error   & order \\ \hline
1/20      &  1.02e-5    &           & 2.58e-5   &        & 1.12e-5     &         & 2.87e-5     &      \\
1/40      &  3.23e-7    & 4.98      & 7.82e-7   &  5.04  & 3.56e-7     & 4.98    & 9.05e-7     & 4.99\\
1/80      &  1.02e-8    & 4.98      & 2.43e-8   &  5.01  & 1.13e-8     & 4.98    & 2.78e-8     & 5.02\\
1/160     &  3.28e-10   & 4.96      & 7.18e-10  &  5.08  & 3.64e-10    & 4.96    & 8.39e-10    & 5.05\\
1/320     &  1.09e-11   & 4.91      & 2.12e-11  &  5.08  & 1.21e-11    & 4.91    & 2.49e-11    & 5.07\\
\bottomrule
\end{tabular}\label{tb52}
\end{table}

\begin{table}[!htbp]\centering
\caption{Error table of our method and that in \cite{TWSN2012} for the linear equation \eqref{51} with solution \eqref{53}.
Here we use the seventh-order  scheme with ideal weights and the SSP(3,3) scheme with $\dt = \dx^{7/3}$.}
\begin{tabular}{lllll | llll}
\toprule
         &  our method  &      &      &      &   method in \cite{TWSN2012} &   &  \\\hline
$\Delta x$ & $L^1$ error & order   & $L^{\infty}$ error   & order & $L^1$ error & order  & $L^{\infty}$ error   & order \\ \hline
1/10      &  1.10e-5    &           & 3.56e-5   &        & 1.17e-5     &         & 3.79e-5     &      \\
1/20      &  9.33e-8    & 6.88      & 3.38e-7   &  6.72  & 9.53e-8     & 6.94    & 3.53e-7     & 6.75 \\
1/40      &  9.54e-10   & 6.95      & 3.77e-9   &  6.49  & 7.64e-10    & 6.96    & 3.88e-9     & 6.51\\
1/80      &  7.67e-12   & 6.62      & 5.97e-11  &  5.98  & 7.81e-12    & 6.61    & 6.01e-11    & 6.01\\
\bottomrule
\end{tabular}\label{tb53}
\end{table}

Fixing $\dx = 1/80$, we  plot in Fig.~\ref{51} the $L^1$ error with respect to the CFL number. We see from the figure that,
with our boundary treatment method,
the errors for all the RK schemes  remain small until a critical CFL number is reached, after which the error grows sharply.
As expected, the RK schemes with negative coefficients give better stability than those without negative coefficients.
\begin{figure}[!ht]
\centering
\includegraphics[width=0.5\textwidth]{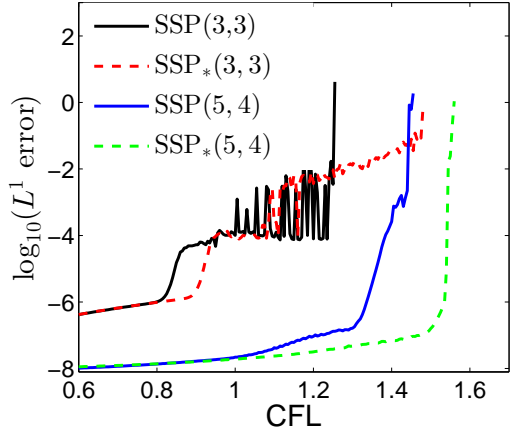}
\caption{$L^1$ errors as a function of the CFL number for the linear system \eqref{51} with solution \eqref{53}.}
\label{fig51}
\end{figure}

\begin{figure}[!ht]
\centering
\includegraphics[width=0.45\textwidth]{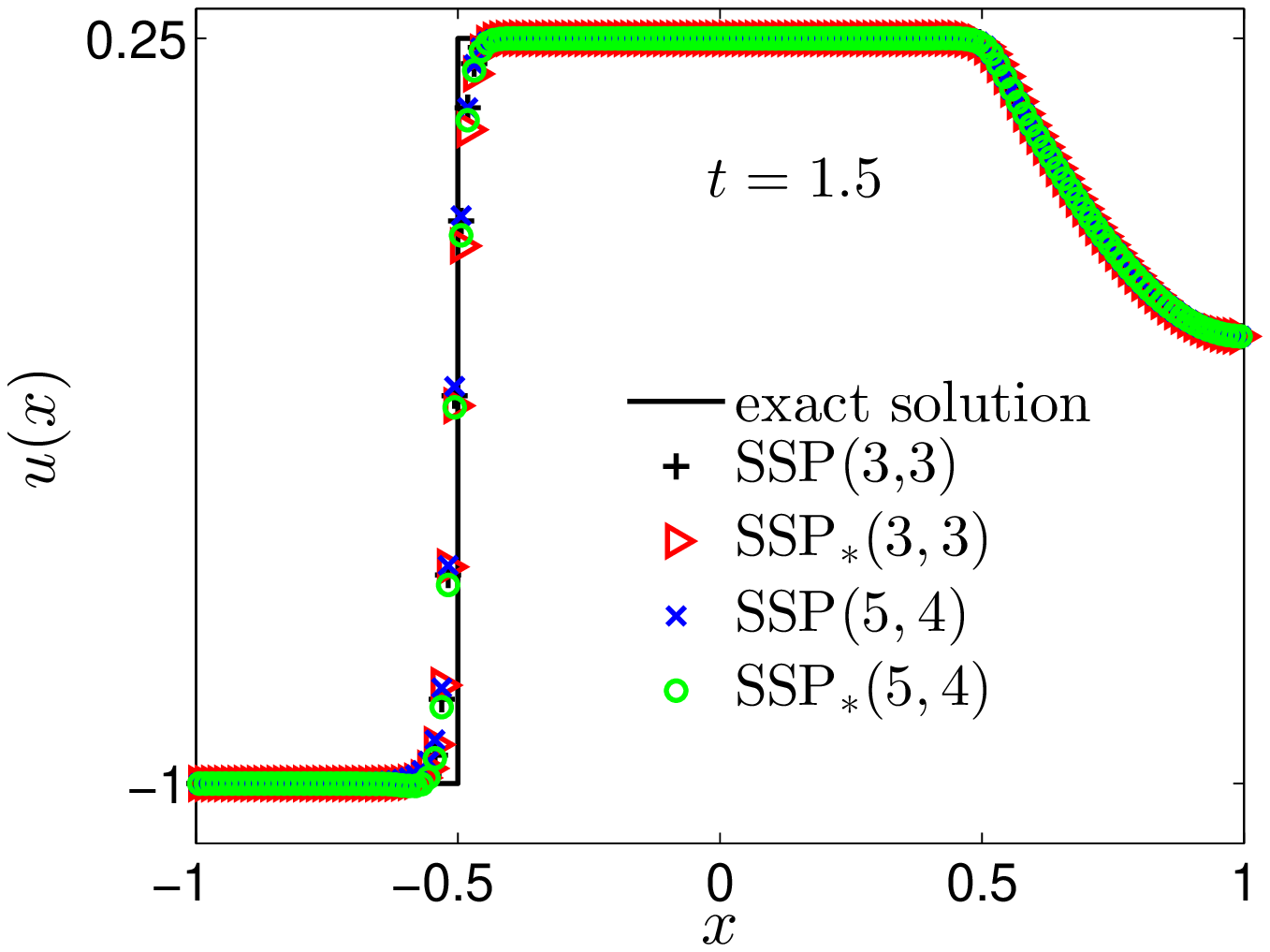}
\includegraphics[width=0.45\textwidth]{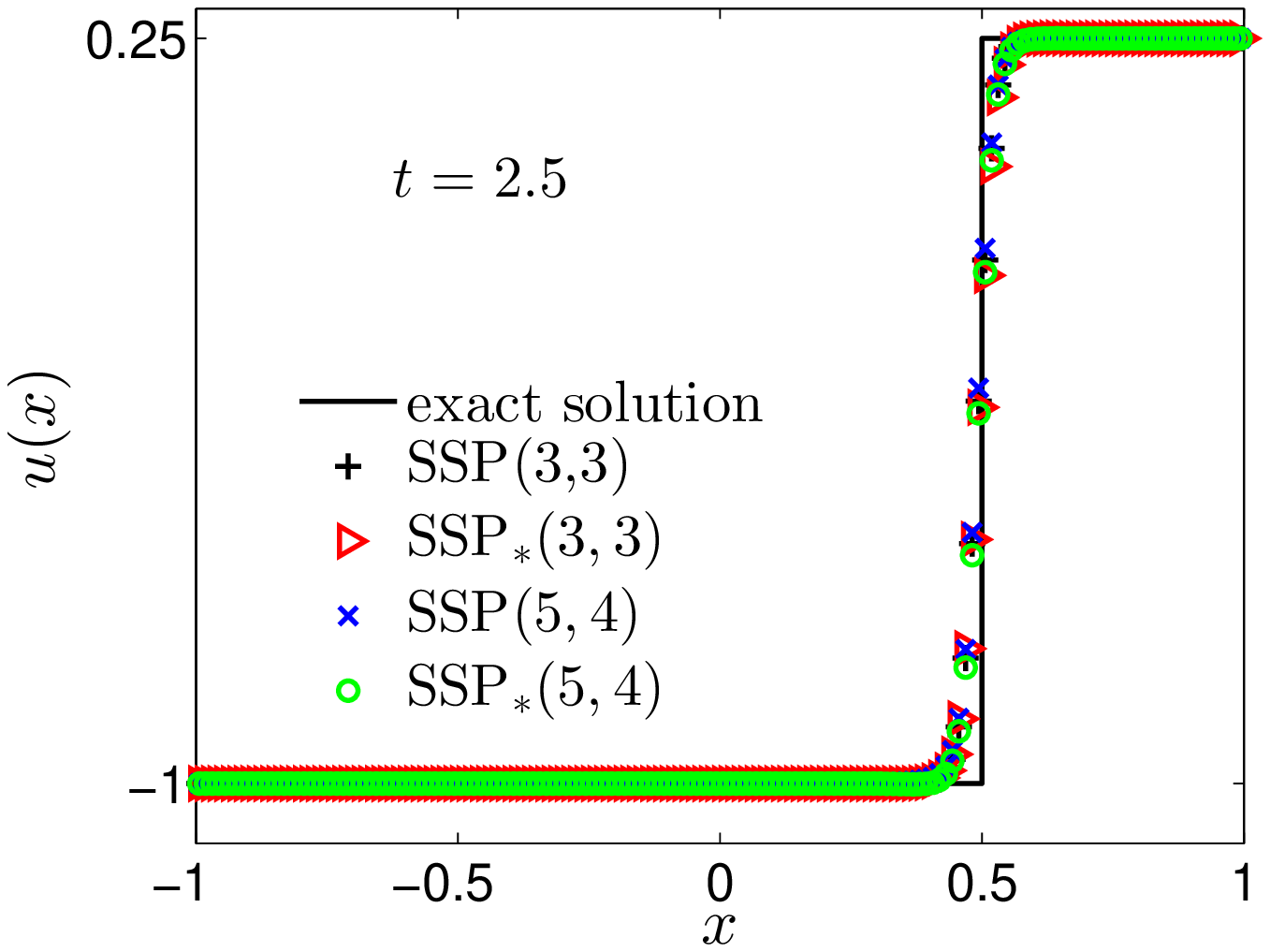}
\caption{Computational results for the linear scalar equation \eqref{51} with boundary condition \eqref{54}. Here we take $\Delta x = 1/80$.}
\label{fig52}
\end{figure}

Next we take \cite{TS2010}
\begin{equation}\label{54}
g(t) =
\left\{
\begin{array}{ll}
 0.25 & t\leq 1,\\
 -1 & t>1
\end{array}
\right.
\end{equation}
for which the exact solution is
\begin{equation}
u(t,x) =
\left\{
\begin{array}{ll}
 -1 & x < t-2, \\
 0.25 & t-2 \leq x < t-1,\\
 0.25 + 0.5 \sin[\pi(x-t)] & x \geq t-1.
\end{array}
\right.
\end{equation}
For this definition of $g(t)$, the exact solution has a slope discontinuity for $t < 2$.  A second discontinuity (this time in the solution) enters the computational domain  from the inflow boundary at $t=1$
and persists until $t=3$.
It can be observed from Fig.~\ref{fig52} (left) that both types of discontinuities are well captured by our method for different RK methods.
These results are comparable to those of the method in \cite{TS2010}.
Note that the computational results remain good after the slope discontinuity passes through the right boundary; see Fig.~\ref{fig52} (right).

%
%
%
%

\subsection{1D Burgers equation}

The second problem is the Burgers equation
\begin{equation}\label{56}
u_t + \left( \frac{1}{2}u^2\right)_x = 0,    \quad  t>0
\end{equation}
on the domain $x \in [-\frac{1}{2},\frac{3}{2}]$.  For $t<1$, this equation has the nonsmooth solution
\begin{equation}\label{57}
u(t,x) =
\left\{
\begin{array}{cl}
 1 & x < t, \\
 \frac{1-x}{1-t} & t \leq x < 2-t,\\
 -1 & x \geq 2-t.
\end{array}
\right.
\end{equation}
In our computations,  initial and boundary data are assigned according to this exact solution.
Note that boundary conditions should be prescribed at both sides here since $u(t,-\frac{1}{2})>0$ and $u(t,\frac{3}{2})<0$.

\begin{figure}[!ht]
\centering
\includegraphics[width=0.45\textwidth]{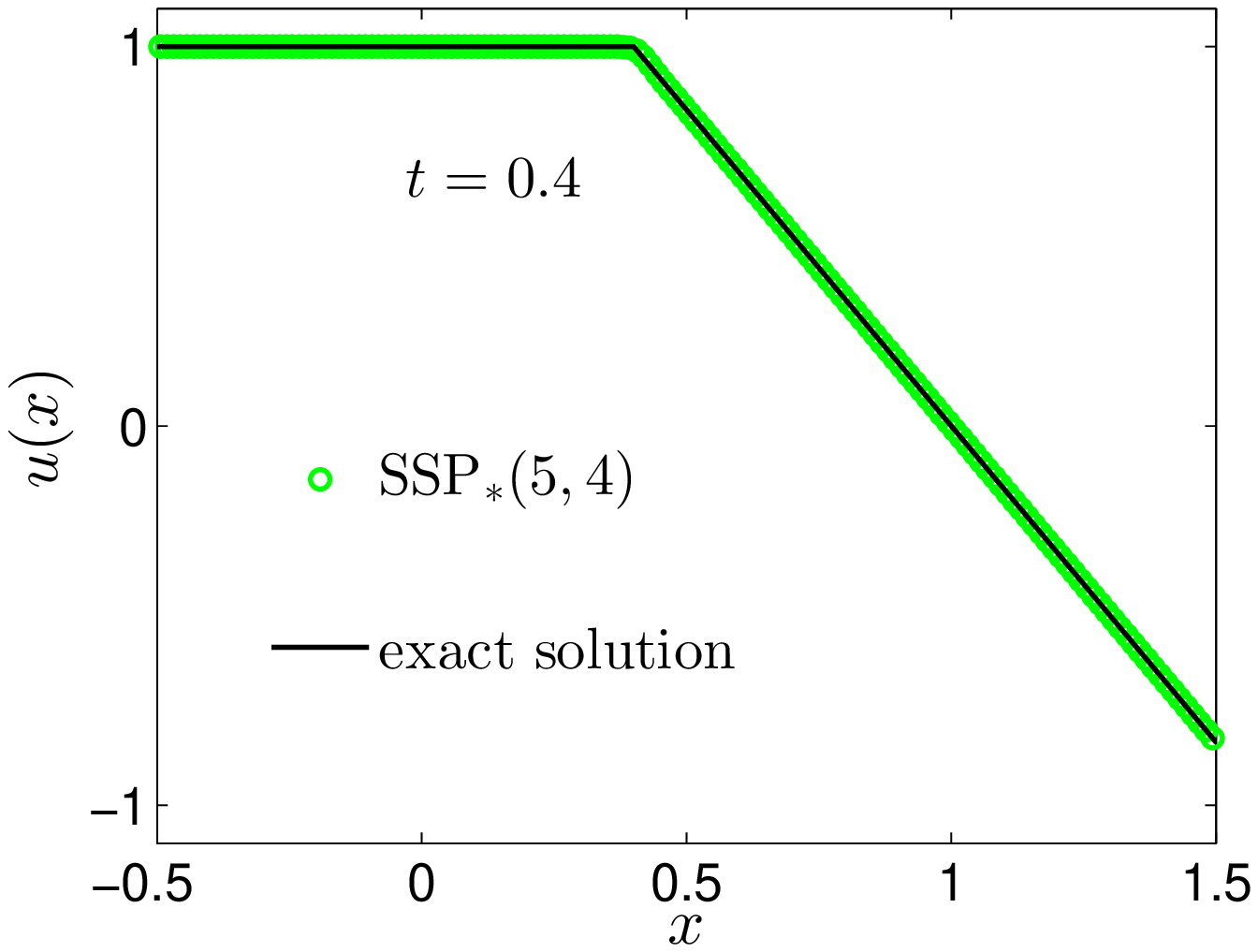}
\includegraphics[width=0.45\textwidth]{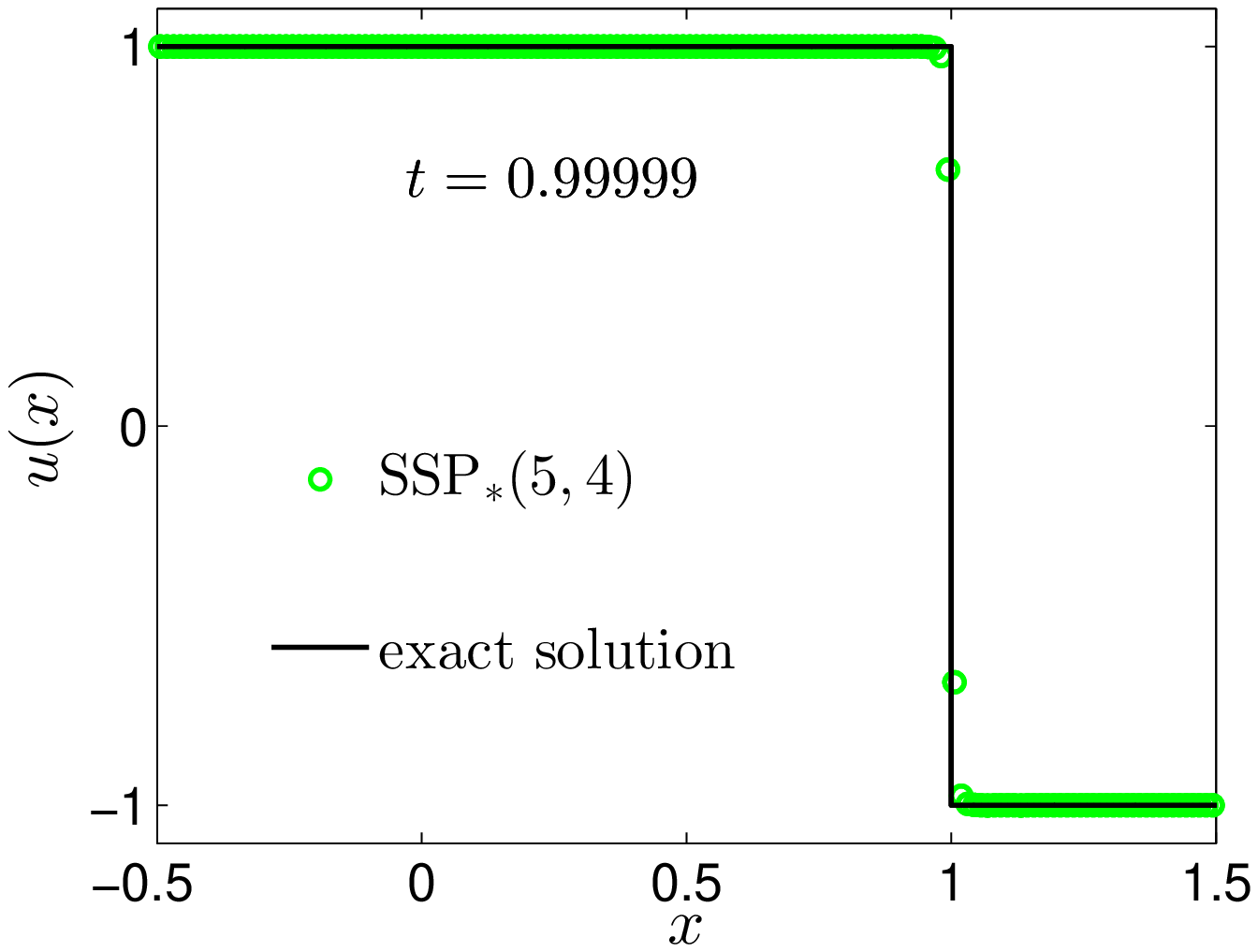}
\caption{Computational results for the Burgers equation \eqref{56} with exact solution \eqref{57}. Here we take $\Delta x = 1/80$.}
\label{fig53}
\end{figure}

With our boundary treatment method used at both boundaries, the computational results of the four RK schemes all agree well with the exact solution.
To illustrate, solutions of the $\mbox{SSP}_{*}(5,4)$  scheme with CFL=0.6  are presented in Fig.~\ref{fig53} for $t=0.4$ and $t=0.99999$. It can be seen that good results are obtained even after the discontinuity has entered from the right-hand boundary.

Next, we compare the errors at $t=0.4$ of the four RK schemes as the CFL number is varied.  For this discontinuous, nonlinear problem, the plots in Fig.~\ref{54} show that the RK schemes with negative coefficients still give better stability than those without negative coefficients,
though the advantage is not as pronounced as that for the smooth solution of the linear equation.
\begin{figure}[!ht]
\centering
\includegraphics[width=0.5\textwidth]{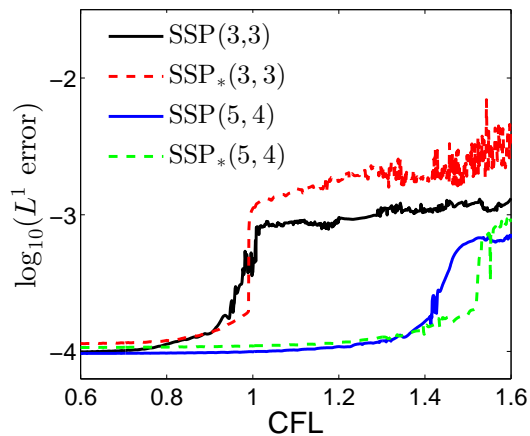}
\caption{$L^1$ errors as a function of the CFL number for the Burgers equation \eqref{56} with nonsmooth solution \eqref{57}.}
\label{fig54}
\end{figure}

\subsection{1D Euler equations}

To conclude our experiments in 1D, we consider the 1D Euler equations
\begin{equation}\label{58}
\left(
\begin{array}{c}
\rho \\
\rho u \\
E
\end{array}
\right)_t
+
\left(
\begin{array}{c}
\rho u \\
\rho u^2 + p \\
u(E+p)
\end{array}
\right)_x = 0,
\end{equation}
where $\rho, u, p$ and $E$ are the fluid density, velocity, pressure and total energy, respectively.
The equation of state has the form
$$
E = \frac{p}{\gamma-1} + \frac{1}{2}\rho u^2,
$$
where $\gamma=1.4$ is the specific heat ratio. The eigenvalues of the Jacobian matrix of the flux are $u-c, u$ and $u+c$ with $c = \sqrt{\gamma p / \rho}$ giving the sound speed.

\subsubsection{Accuracy test}

We first test the accuracy of our method for a smooth solution on the domain $[-\pi, \pi]$ \cite{TS2010}:
\begin{equation}\label{59}
\begin{split}
& \rho(t,x) = 1+0.2\sin(x-t), \\
& u(t,x) = 1, \\
& p(t,x) = 2.
\end{split}
\end{equation}
The initial condition is obtained by setting $t=0$ in the above solution. It can be directly verified that
$u-c<0$, $u>0$ and $u+c>0$ at both boundaries. Thus, two boundary conditions are needed at $x=-\pi$ and one is required at $x=\pi$. As in \cite{TS2010}, we prescribe $\rho$ and $u$ at the left boundary:
\begin{equation*}
\begin{split}
& \rho(t,-\pi) = 1+0.2\sin(t), \\
& u(t,-\pi) = 1,
\end{split}
\end{equation*}
and $\rho$ is given at the right boundary:
\begin{equation*}
\rho(t,\pi) = 1+0.2\sin(t).
\end{equation*}
The errors and convergence orders of different RK schemes using a fixed CFL number of 0.6 are listed in Table \ref{tb54} for this example.
We see that for the third-order RK schemes, the convergence order is much larger than three for coarse meshes and that third-order
convergence is achieved by taking very small $\dx$-values. Similar to the linear problem, the convergence order is about five for the
fourth-order RK schemes.  As before, we conjecture that this is due to  the error of the space discretization dominating over that of the time discretization.
\begin{table}[!htbp]\centering
\caption{Error table for the Euler equations \eqref{58} with solution \eqref{59} at $t=2$. }
\begin{tabular}{lllll | llll}
\toprule
         &  SSP(3,3)  &      &      &      &   $\mbox{SSP}_{*}(3,3)$ &   &  \\\hline
$\Delta x$ & $L^1$ error & order   & $L^{\infty}$ error   & order & $L^1$ error & order  & $L^{\infty}$ error   & order \\ \hline
$\pi/20$      &  5.56e-6    &           & 1.62e-5   &        & 8.33e-6     &         & 2.31e-5     &      \\
$\pi/40$      &  1.89e-7    & 4.88      & 5.61e-7   &  4.85  & 2.75e-7     & 4.92    & 7.88e-7     & 4.87\\
$\pi/80$      &  8.85e-9    & 4.42      & 2.38e-8   &  4.56  & 1.15e-8     & 4.58    & 3.06e-8     & 4.69\\
$\pi/160$     &  6.42e-10   & 3.79      & 1.53e-9   &  3.96  & 7.26e-10    & 3.99    & 1.73e-9     & 4.14\\
$\pi/320$     &  6.57e-11   & 3.29      & 1.51e-10  &  3.34  & 6.82e-11    & 3.41    & 1.56e-10    & 3.47\\
$\pi/640$     &  7.75e-12   & 3.08      & 1.81e-11  &  3.06  & 7.82e-12    & 3.12    & 1.81e-11    & 3.11\\
\toprule
         &  SSP(5,4)  &      &      &      &   $\mbox{SSP}_{*}(5,4)$ &   &  \\\hline
$\Delta x$ & $L^1$ error & order   & $L^{\infty}$ error   & order & $L^1$ error & order  & $L^{\infty}$ error   & order \\ \hline
$\pi/20$      &  5.33e-6    &           & 1.57e-5   &        & 7.02e-6     &         & 2.00e-5     &      \\
$\pi/40$      &  1.59e-7    & 5.07      & 4.91e-7   &  5.00  & 2.11e-7     & 5.06    & 6.31e-7     & 4.99\\
$\pi/80$      &  5.00e-9    & 4.99      & 1.48e-8   &  5.05  & 6.61e-9     & 4.99    & 1.90e-8     & 5.05\\
$\pi/160$     &  1.56e-10   & 5.00      & 4.13e-10  &  5.16  & 2.07e-10    & 5.00    & 5.32e-10    & 5.16\\
$\pi/320$     &  5.01e-12   & 4.96      & 1.26e-11  &  5.03  & 6.33e-12    & 5.03    & 1.53e-11    & 5.12\\
\bottomrule
\end{tabular}\label{tb54}
\end{table}

\subsubsection{Comparison with the method in \cite{TS2010}}

We now verify our analysis  of Section~4  where  it was shown that order reduction arises when the method in \cite{TS2010} is applied to the nonlinear equations.
To this end, we apply the SSP(3,3) RK method and the seventh-order WENO scheme to the 1D Euler equations with smooth solution \eqref{59}.
As in the computation for the  linear scalar equation,  the ideal WENO weights are used.
Applying a fixed $\dt = \dx^{7/3}$, we compare in Table \ref{tb55} the errors of our method and that in \cite{TS2010}. It is clear that the designed seventh-order convergence is obtained for our method while the convergence order of the method in \cite{TS2010} is only about $\frac{17}{3}(\approx 5.67)$.  This agrees with our analysis:
 our boundary treatment method maintains the accuracy of the interior spatial and time discretizations while
 the method in \cite{TS2010} introduces order reduction.

\begin{table}[!htbp]\centering
\caption{Error table of our method and that in \cite{TWSN2012} for the Euler equations \eqref{58} with solution \eqref{59} at $t=2$.
Here we use the seventh-order WENO scheme with ideal weights and the SSP(3,3) scheme with $\dt = \dx^{7/3}$.}
\begin{tabular}{lllll | llll}
\toprule
         &  our method  &      &      &      &   method in \cite{TWSN2012} &   &  \\\hline
$\Delta x$ & $L^1$ error & order   & $L^{\infty}$ error   & order & $L^1$ error & order  & $L^{\infty}$ error   & order \\ \hline
$\pi/10$      &  5.09e-6    &           & 1.31e-5   &        & 1.41e-5     &         & 4.52e-5     &      \\
$\pi/20$      &  4.36e-8    & 6.87      & 1.23e-7   &  6.73  & 2.16e-7     & 6.03    & 8.04e-7     & 5.81 \\
$\pi/40$      &  3.34e-10   & 7.03      & 1.00e-9   &  6.94  & 4.06e-9     & 5.73    & 1.49e-8     & 5.75\\
$\pi/80$      &  3.54e-12   & 7.04      & 7.37e-12  &  7.08  & 7.62e-11    & 5.73    & 2.88e-10    & 5.69\\
\bottomrule
\end{tabular}\label{tb55}
\end{table}

\subsubsection{Stability test}

\begin{figure}[!ht]
\centering
\includegraphics[width=0.5\textwidth]{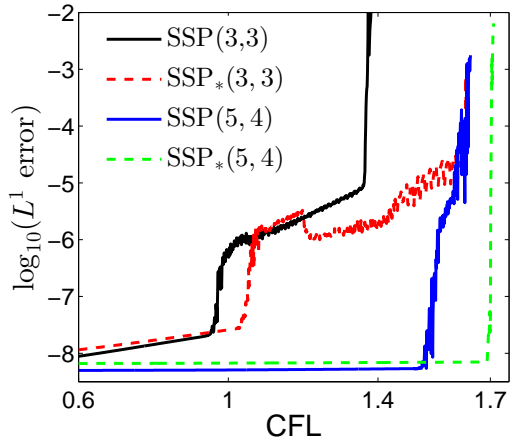}
\caption{$L^1$ errors as a function of the CFL number for the Euler equations \eqref{58} with solution \eqref{59}.}
\label{fig55}
\end{figure}

We conduct another test of the stability of our boundary treatment method for the Euler equations with smooth solution~\eqref{59} by fixing $\dx = 1/80$ and varying the CFL number.
The corresponding $L^1$ errors are plotted in Fig.~\ref{fig55}.
Similar to the case of the linear equation, we see that the RK schemes with negative coefficients give better stability than those without negative coefficients.
Furthermore, with our boundary treatment method,
the error for all RK schemes remains small until a critical CFL number is reached, after which the error grows sharply.

\subsubsection{Performance for shock interactions}

To conclude, we test our method for the Euler equations on a problem with shocks that arise through the interaction of two blast waves \cite{WC2000,TS2010}.
In this problem, the initial data are given by
\begin{equation}
(\rho, u, p)(0,x) =
\left\{
\begin{array}{ll}
(1,0,10^3)    & \quad 0<x<0.1,\\
(1,0,10^{-2}) & \quad 0.1<x<0.9,\\
(1,0,10^2)    & \quad 0.9<x<1
\end{array}
\right.
\end{equation}
and the boundaries are solid walls, namely, $u=0$ at $x=0$ and $x=1$. As time evolves, there are multiple reflections of shocks and rarefactions off the walls. The shocks and rarefactions also interact with each other and with contact discontinuities. Similar to \cite{TS2010}, we use second-order Taylor expansion at the boundary.  The WENO type extrapolation in \cite{TS2010} is employed to avoid oscillations.

\begin{figure}[!ht]
\centering
\includegraphics[width=0.45\textwidth]{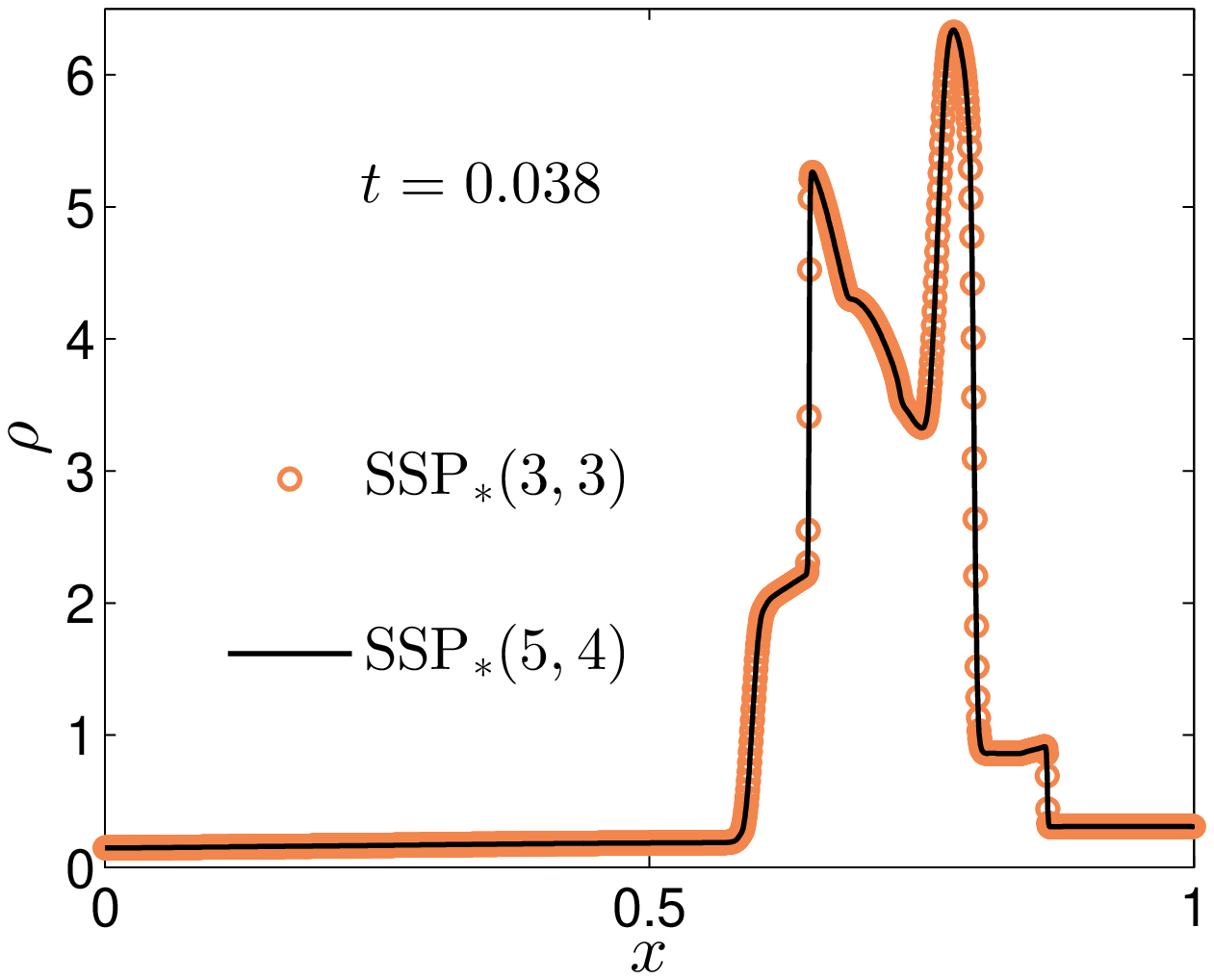}
\includegraphics[width=0.45\textwidth]{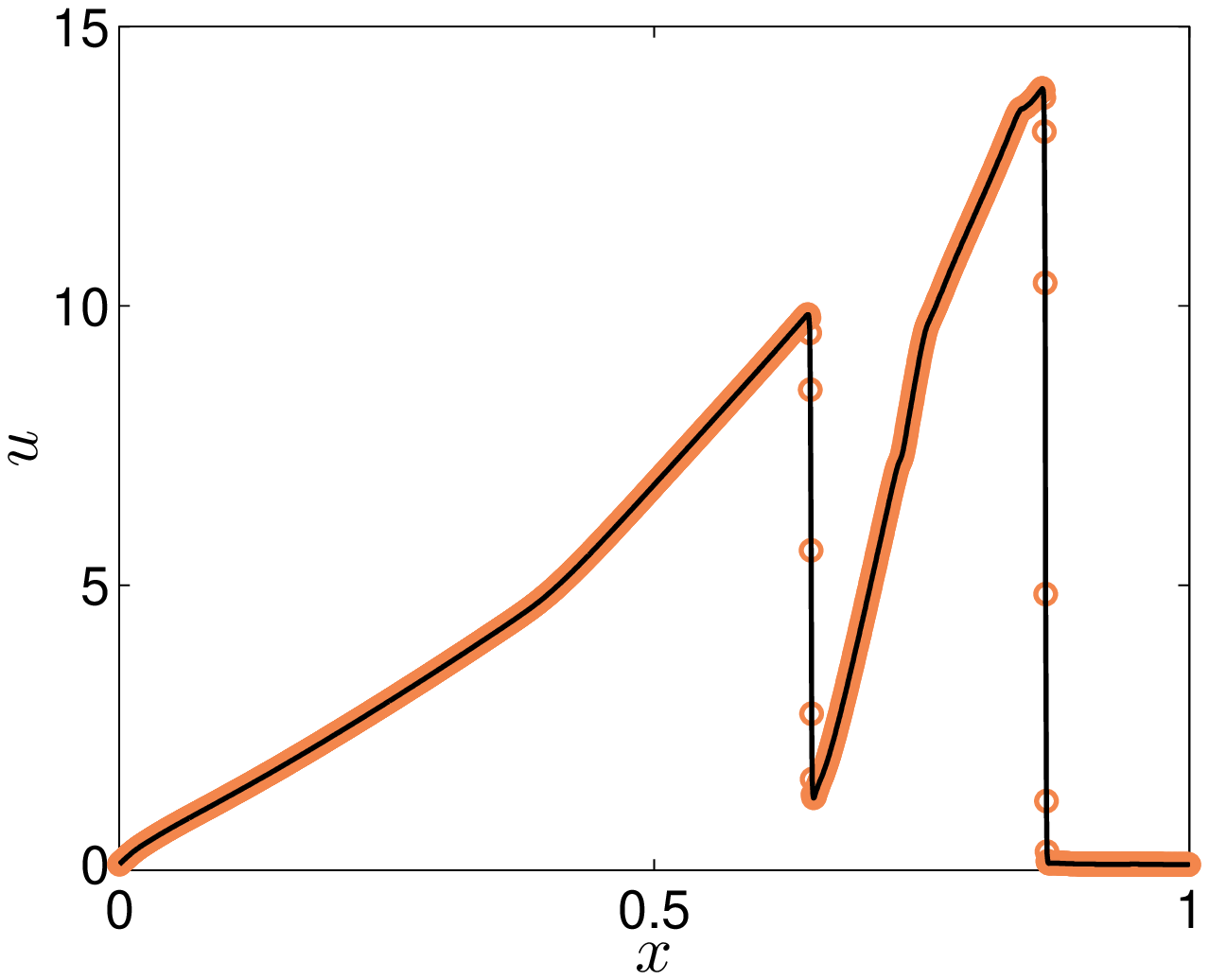}
\includegraphics[width=0.45\textwidth]{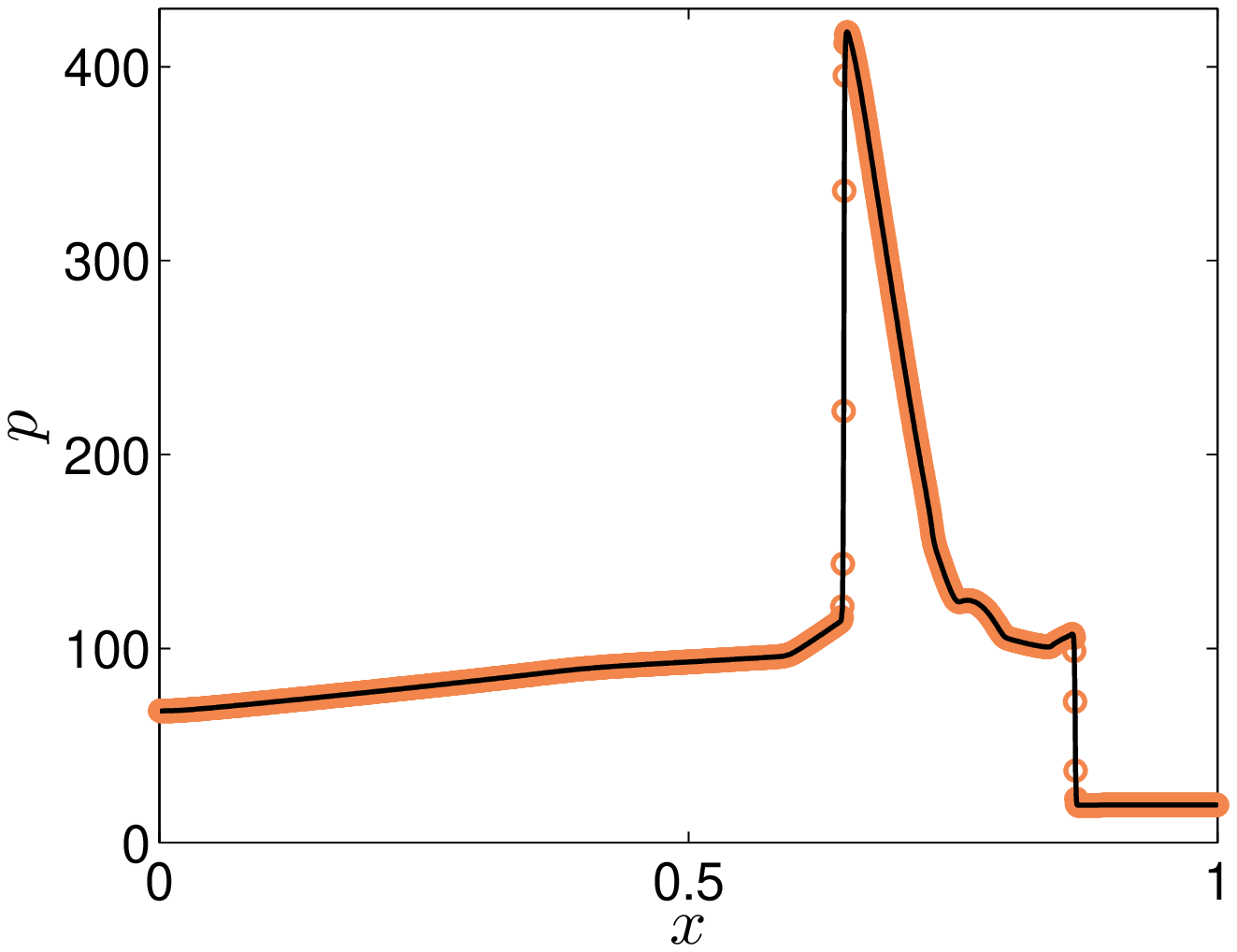}
\caption{Computational results of the blast wave problem for the Euler equations \eqref{56} with $\dx = 1/1600$.}
\label{fig56}
\end{figure}

Fig.~\ref{fig56} gives a plot of
our numerical solutions at $t=0.038$.
We see that the results of the $\mbox{SSP}_{*}(3,3)$ and $\mbox{SSP}_{*}(5,4)$ schemes coincide with each other.  They also agree well with those of the method in \cite{TS2010}. These good results demonstrate the ability of our boundary treatment method to compute solutions
to problems with complex shock interactions.

\subsection{2D Euler equations}

Though our method is illustrated only for 1D equations in Section \ref{sec3}, it can be directly extended to 2D problems \cite{ZH2019}. In the following we examine our method for the 2D Euler equations
\begin{equation}\label{eq:euler}
\partial_t U + \partial_x F(U) + \partial_y G(U) = 0, \quad  (x,y) \in \Omega
\end{equation}
with
\begin{equation*}
\begin{split}
& U = (\rho, \rho u, \rho v, E)^T, \\
& F(U) = ( \rho u, \rho u^2 + p, \rho u v, (E+p)u)^T, \\
& G(U) = (\rho v, \rho u v, \rho v^2 + p, (E+p)v)^T
\end{split}
\end{equation*}
and equation of state $E=\frac{1}{2}\rho(u^2 + v^2) + \frac{p}{\gamma-1}$.
Here $\rho$ is the density, $u$ and $v$ are the velocities in $x$ and $y$ directions, $E$ is the total energy, $p$ is the pressure, and $\gamma=1.4$ is the specific heat ratio. Specifically, we simulate the vortex problem with analytical solutions to the above 2D Euler equations as in \cite{TS2010}. In this problem, the mean flow $\rho=u=v=1$ is initially imposed with an isentropic vortex perturbation centered at $(x_0,y_0)$ in $(u,v)$ and in the temperature $T=p/\rho$, no perturbation in the entropy $S=p/\rho^\gamma$:
\begin{equation*}
\begin{split}
& (\delta u, \delta v) = \frac{\epsilon}{2\pi}e^{0.5(1-r^2)}(-\bar y, \bar x), \\
& \delta T = -\frac{(\gamma-1)\epsilon^2}{8\gamma\pi^2}e^{(1-r^2)}, \\
& \delta S = 0.
\end{split}
\end{equation*}
Here $(\bar x, \bar y) = (x-x_0, y-y_0)$, $r^2=\bar x^2 + \bar y^2$ and $\epsilon$ is the vortex strength. The exact solution $U(t,x,y)$ of this problem is a passive convection of the vortex with the mean velocity, \ie, $U(t,x,y) = U(0,x-ut,y-vt)$. As in \cite{TS2010}, we set the spatial domain $\Omega = [-0.5, 1]^2$ and $(x_0,y_0)=(0,0)$, and take boundary conditions from the exact solution whenever needed.

In the computation, we use the fifth-order WENO scheme \cite{SO1988} for the spatial discretization. The tangential derivative $\partial_y U^n$ at the boundary is given by the exact solution and $\partial_y U^{(i)}, \partial_{xy} U^n, \partial_{xy} U^{(i)}$ are computed by numerical differentiation. We take $\epsilon = 1$ so that the signs of the eigenvalues at the boundary, and thereby the number of boundary conditions, do not change before the terminal time $t=1$. The recently proposed method in \cite{LSTZ2020} may be used to deal with the situation where the signs of the eigenvalues vary at the boundary. In this case, the errors are given in Table \ref{tb56}. It can be seen that the convergence orders of the computational results are generally larger than those of the four RK schemes. This may be due to the error of the spatial discretization dominating over that of the time discretization. These results verify the effectiveness of our method for 2D problems.

\begin{table}[!htbp]\centering
\caption{Error table for the vortex problem of the 2D Euler equations at $t=1$ with CFL=0.6. Fifth-order WENO scheme is used for the spatial discretization.}
\begin{tabular}{lllll | llll}
\toprule
         &  SSP(3,3)  &      &      &      &   $\mbox{SSP}_{*}(3,3)$ &   &  \\\hline
$\Delta x$ & $L^1$ error & order   & $L^{\infty}$ error   & order & $L^1$ error & order  & $L^{\infty}$ error   & order \\ \hline
$1.5/20$      &  6.17e-7    &           & 6.09e-6   &        & 1.00e-6     &         & 1.62e-5     &      \\
$1.5/40$      &  2.11e-8    & 4.87      & 3.94e-7   &  3.95  & 3.35e-8     & 4.91    & 1.69e-6     & 3.25\\
$1.5/80$      &  7.83e-10   & 4.75      & 1.62e-8   &  4.61  & 1.61e-9     & 4.38    & 3.65e-7     & 2.21\\
$1.5/160$     &  4.83e-11   & 4.02      & 1.76e-10  &  6.52  & 5.71e-11    & 4.82    & 8.21e-9     & 5.47\\
\toprule
         &  SSP(5,4)  &      &      &      &   $\mbox{SSP}_{*}(5,4)$ &   &  \\\hline
$\Delta x$ & $L^1$ error & order   & $L^{\infty}$ error   & order & $L^1$ error & order  & $L^{\infty}$ error   & order \\ \hline
$1.5/20$      &  5.89e-7    &           & 4.71e-6   &        & 7.58e-7     &         & 8.38e-6     &      \\
$1.5/40$      &  1.93e-8    & 4.93      & 2.62e-7   &  4.17  & 2.42e-8     & 4.97    & 6.37e-7     & 3.72\\
$1.5/80$      &  6.04e-10   & 5.00      & 7.52e-9   &  5.12  & 7.55e-10    & 5.00    & 4.87e-8     & 3.71\\
$1.5/160$     &  2.82e-11   & 4.42      & 1.99e-10  &  5.24  & 3.06e-11    & 4.62    & 2.75e-10    & 7.47\\
\bottomrule
\end{tabular}\label{tb56}
\end{table}

\section{Conclusions and remarks}

In this paper, we examine the finite difference boundary treatment method proposed in our previous work \cite{ZH2019} for the case of explicit high-order SSP RK schemes applied to hyperbolic conservation laws.
Our spatial discretization uses a WENO scheme on a Cartesian mesh where the corresponding downwind scheme is applied in the case of negative coefficients. We show that our method applies to general SSP RK schemes with/without negative coefficients.  Our method preserves the accuracy of the RK schemes and has good stability. These nice properties are verified on linear and nonlinear problems using third- and fourth-order RK schemes with and without negative coefficients.
In addition, when boundary conditions are present and our  boundary treatment method is used, the SSP RK schemes with negative coefficients still allow for larger time steps than those with all nonnegative coefficients. In this regard, the boundary treatment method in \cite{ZH2019} is an effective supplement to SSP RK schemes with/without negative coefficients for initial-boundary value problems of conservation laws.

It should be remarked that the key point of our method is to use the RK scheme itself at the boundary, instead of imposing intermediate boundary conditions as in \cite{TS2010,TWSN2012}.
The boundary conditions in those previous works~\cite{TS2010,TWSN2012} are derived only for RK schemes up to third order while our method applies to RK schemes of arbitrary order. In addition, we show the equivalence of our method and that in \cite{TS2010,TWSN2012} for linear problems and the case of the third-order SSP RK scheme.  On the other hand, for nonlinear equations, there exist difference terms of $\order(\dt^2\dx)$ between the two methods, which result in order reduction for the method in \cite{TS2010,TWSN2012}. These results have been numerically verified as well.

\section*{Acknowledgements}
This work was supported by the National Foreign Experts Project of China (No. G20190001349) and the National Natural Science Foundation of China
(No. 11801030, No. 11861131004).

\section*{Appendix}
This appendix gives the computation of $U_j^{(2)}$ at ghost points $j=-1,-2,-3$. As for $U_j^{(1)}$, $U_j^{(2)}$ is approximated with the fifth-order Taylor expansion at the boundary point $x_b=0$:
\begin{equation}\label{a1}
U_j^{(2)} = \sum_{k=0}^{4} \frac{x_j^k}{k!} U^{(2),(k)},  \quad j=-1,-2,-3,
\end{equation}
where $U^{(2),(k)}$ denotes a $(5-k)$-th order approximation of the spatial derivative $\frac{\partial^k U^{(2)}}{\partial x^k}\big|_{x=0}$.
Next, we compute $U^{(2),(k)}$ for $k=0,1,2,3,4$.
To this end, we apply the second stage of the RK solver \eqref{24b} for $U^{(2)}$ at the boundary point $x_b=0$:
\begin{equation}\label{a2}
U^{(2)}(x_b) = \alpha_{20} U^{n}(x_b) + \alpha_{21} U^{(1)}(x_b)
- \beta_{20}  \Delta t \partial_x F(U^n(x_b))
- \beta_{21}  \Delta t \partial_x F(U^{(1)}(x_b)).
\end{equation}
Notice that
\begin{equation*}
\begin{split}
& U^{n}(x_b) = U^{n,(0)} + \order(\Delta x^5), \\
& U^{(1)}(x_b) = U^{(1),(0)} + \order(\Delta x^5), \\
& \partial_x F(U^n(x_b)) = F_U(U^n(x_b)) \partial_x U^n(x_b) = F_U(U^{n,(0)}) U^{n,(1)} + \order(\Delta x^4),\\
& \partial_x F(U^{(1)}(x_b)) = F_U(U^{(1)}(x_b)) \partial_x U^{(1)}(x_b) = F_U(U^{(1),(0)}) U^{(1),(1)} + \order(\Delta x^4).
\end{split}
\end{equation*}
Substituting these into \eqref{a2}, we obtain an approximation of  $U^{(2)}(x_b)$ and denote it by $U^{(2),(0)}$. Observe that the error between $U^{(2),(0)}$ and  $U^{(2)}(x_b)$ defined by \eqref{a2} is $\order(\Delta x^5)$.

Taking derivatives with respect to $x$ on both sides of \eqref{a2} yields
\begin{equation}\label{a3}
\frac{\partial U^{(2)}}{\partial x}\big|_{x= x_b}
= \alpha_{20} \partial_x U^{n} (x_b) + \alpha_{21} \partial_x U^{(1)} (x_b)
 -\beta_{20} \Delta t \partial_{xx} F(U^n(x_b))
 -\beta_{21} \Delta t \partial_{xx} F(U^{(1)}(x_b)).
\end{equation}
With approximations
\begin{equation*}
\begin{split}
& \partial_x U^{n} (x_b) = U^{n,(1)} + \order(\Delta x^4), \\
& \partial_x U^{(1)} (x_b) = U^{(1),(1)} + \order(\Delta x^4), \\
& \partial_{xx} F(U^n(x_b)) = F_{UU}(U^n(x_b)) \partial_x U^n(x_b) \partial_x U^n(x_b) + F_U(U^n(x_b)) \partial_{xx} U^n(x_b) \\
& \hspace{2.45cm}= F_{UU}(U^{n,(0)}) U^{n,(1)} U^{n,(1)} + F_U(U^{n,(0)}) U^{n,(2)} +  \order(\Delta x^3),\\
& \partial_{xx} F(U^{(1)}(x_b)) = F_{UU}(U^{(1)}(x_b)) \partial_x U^{(1)}(x_b) \partial_x U^{(1)}(x_b) + F_U(U^{(1)}(x_b)) \partial_{xx} U^{(1)}(x_b)\\
& \hspace{2.65cm}= F_{UU}(U^{(1),(0)}) U^{(1),(1)} U^{(1),(1)} + F_U(U^{(1),(0)}) U^{(1),(2)} +  \order(\Delta x^3),
\end{split}
\end{equation*}
$\frac{\partial U^{(2)}}{\partial x}\big|_{x= x_b} $ can be computed with \eqref{a3} and the resulting solution  $U^{(2),(1)}$ is a fourth-order approximation of $\frac{\partial U^{(2)}}{\partial x}\big|_{x= x_b} $.

Furthermore, we simply approximate $\frac{\partial ^k U^{(2)}}{\partial x^k}\big|_{x= x_b} $ for $k\geq2$ by using a $(5-k)$-th order extrapolation with $U^{(2)}$ at interior points.
In this way, we obtain $U^{(2),(k)}$ for $k=0,1,2,3,4$ with accuracy of order $(5-k)$. Then, $U_j^{(2)}$ for $j=-1,-2,-3$ can be computed by \eqref{a1} with fifth-order accuracy. Having $U^{(2)}$ at the ghost points, we can then evolve from $U^{(2)}$ to $U^{(3)}$ using the interior difference scheme.


\end{document}